\documentclass[%
 reprint, showkeys,
 amsmath,amssymb,
 aps,
]{revtex4-2}

\usepackage{graphicx}
\usepackage{dcolumn}
\usepackage{bm}
\usepackage{physics}
\usepackage{amsthm}
\usepackage{xcolor}
\usepackage{subcaption}
\usepackage{mathtools}
\usepackage{caption}
\usepackage[caption=false]{subfig}
\newtheorem*{remark}{Remark}

\newcommand{\e}{\ensuremath{\mathrm{e}}}
\newcommand{\letsymb}[1]{#1}
\newcommand{\average}[1]{\left\langle #1\right\rangle}
\begin{document}

\preprint{APS/123-QED}

\title{The Funessian process: Non-Markovian dynamics shaped by the first event}

\author{Bilal Canturk}
 \thanks{corresponding author: bcanturk@gelisim.edu.tr}
\affiliation{Faculty of Architecture and Engineering, Department of Basic Sciences, Istanbul Gelisim University, 34310 Istanbul, Türkiye
}%
\author{G. Baris Bagci}
\affiliation{%
	Department of Physics, Mersin University, 33343 Mersin, Türkiye
}

\author{Onur Pusuluk}
\affiliation{Faculty of Engineering and Natural Sciences, Kadir Has University, 34083 Istanbul, Türkiye.
}

\date{\today}

\begin{abstract}
We construct a continuous-time, positively divisible non-Markovian process with memory of the initial state that satisfies the differential Chapman--Kolmogorov equation. In the stationary state, the correlation function exhibits exponential decay, a behavior typically regarded as characteristic of Markovian dynamics. Nevertheless, the memory is preserved throughout the evolution of the process, manifesting itself in observable statistical quantities. We further demonstrate that mutual information serves as a reliable measure of the non-Markovian character of the process. As an application, we study a random walk driven by the constructed process and show that the memory effect breaks ergodicity and modifies transport properties such as the diffusion coefficient.
\end{abstract}
\maketitle


\section{Introduction}
\par It is commonly assumed that the Chapman–Kolmogorov equation ($\mathrm{CK}$-equation), together with its differential form, characterizes Markovian dynamics. However, this assumption was already questioned by Doob in an early discussion with Feller~\cite{Snell1997}. In particular, Feller constructed examples of non-Markov processes satisfying the $\mathrm{CK}$-equation under specific initial distributions~\cite{Feller1959}. Van Kampen provided a similar example~\cite{Kampen2007}. These constructions rely on the highly restrictive choice of a uniform initial distribution and therefore do not constitute general counterexamples. More recently, Canturk and Breuer~\cite{CanturkBreuer2024} presented a discrete-time, two-dimensional non-Markov process that satisfies the $\mathrm{CK}$-equation without such restrictions. While their construction is theoretically valid, it is based on replicating a finite number of initial time steps, which inherently limits its applicability to more general or continuous-time settings.
\par The characterization of Markov processes, especially in continuous-time settings, remains a fundamental problem with significant practical implications for stochastic processes~\cite{MetzlerKlafter2004,KutnerMasoliver2017, Afek2017}, statistical physics~\cite{BouchaudGeorges1990, MarksteinerZoller1996, WangTokuyama1999, Afek2020, Dechant2012,Dechant2014,Dechant2016}, and open quantum systems~\cite{Vacchini-Smirne-Breuer2011,DeVeGa2017a, Bouganne2020, BaoLi2025}. Despite its central role, there is still no universally accepted criterion that fully captures Markovianity~\cite{,CanturkBreuer2024,Li2018a,Utagi2023a}. This lack of a clear characterization has led to multiple, and sometimes inconsistent, formulations of Markovian dynamics, particularly in the context of open quantum systems \cite{Li2018a}.
\par In the present paper, noting that \textit{classical P-divisibility} extends the $\mathrm{CK}$-equation to stochastic matrices~\cite{CanturkBreuer2024, Vacchini-Smirne-Breuer2011}, as stated below, we address the following crucial questions concerning the characterization of Markov processes:
\begin{itemize}
	\item[1)] Is there any continuous-time \textit{classical P-divisible} non-Markovian process, satisfying differential $\mathrm{CK}$-equation? 
	\item[2)]  If so, is it possible to distinguish such processes from Markov processes by some means of measurement?
\end{itemize}      
We answer both questions affirmatively. In particular, inspired by Ref.~\cite{CanturkBreuer2024}, we first establish a theoretical framework in Section~\ref{theoretical_framwork}, based on which we construct a continuous-time, $2$-dimensional, classical P-divisible non-Markovian process in Section~\ref{construct_the_process} which retains memory of the initial state and satisfies differential  $\mathrm{CK}$-equation. Secondly, after providing its statistical analysis, we demonstrate in Section~\ref{statistical_anlaysis} that mutual information can serve as a reliable measure of non-Markovianity. In Section~\ref{application_random_walk} we apply the constructed process to random walk. Finally, we discuss our result in the concluding Section~\ref{conclusion_section}. 
\section{Theoretical framework for constructing the process}\label{theoretical_framwork}

We consider a general stochastic process $\{X_{t}\}_{t\in\mathrm{T}}$ with sample space $E=\{x_1,x_2,\ldots,x_d\}$ evolving continuously in time. Furthermore, let $\mathbf{q}(t_0)=(q_1,q_2,\ldots,q_d)$ be the initial probability at time $t_{0}$. Then, the probability vector $\mathbf{p}(t)=(p(x_1,t),p(x_2,t),\ldots,p(x_d,t))$ at a later time $t\geq\letsymb{t_0}$ is determined by the equality $\mathbf{p(t)}=\Lambda(t\vert\letsymb{t_0})\mathbf{q}(t_0)$, where $\Lambda(t\vert\letsymb{t_0})$ is the transition matrix defined in terms of the conditional probabilities as $\Lambda(t\vert\letsymb{t_0})_{jk}:=p(x_j,t\vert\letsymb{x_k},t_0)$. Thus, $p(x_j,t;\letsymb{x_k},t_{0})=\Lambda(t\vert\letsymb{t_0})_{jk}q_{k}$. A general transition matrix $\Lambda(t\vert\letsymb{s})$ is characterized  by the following properties for all $t_0\leq\letsymb{s}\leq\letsymb{t}$: 
\begin{enumerate}
	\item [1.] Its elements are non-negative, $\Lambda(t\vert\letsymb{s})_{jk}\geq 0$,
	\item [2.] $\sum_{j=1}^{d}  \Lambda(t\vert\letsymb{s})_{jk}=1$ for all $k$,
	\item [3.] $\lim_{t\rightarrow\letsymb{s}}\Lambda(t\vert\letsymb{s})=I$, where $I$ is the identity matrix,
	\item [4.] $p_2(x_{j},t;x_{k},s)=\Lambda(t \vert\letsymb{s})_{jk}p(x_{k},s)$ for all $s \leq\letsymb{t}$.
\end{enumerate}  
Dropping the fourth property, the notion of a transition matrix generalizes to that of a stochastic matrix, which is no longer subject to this constraint.
\par  The joint probability of order three for $t_{0}\leq\letsymb{s}\leq\letsymb{t}$ takes the form 
\begin{equation}\label{joint_prob_order_3_eq1}
		p(x_j,t;x_k,s;x_{l},t_{0})=p(x_{j},t\vert\letsymb{x}_{k},s;x_{l},t_{0})\Lambda(s\vert\letsymb{t}_{0})_{kl}q_{l}
\end{equation}
 We recall the \textit{one-point memory transition matrices}  introduced in Ref.~\cite{CanturkBreuer2024},
\[Q_{1}^{(l)}(t\vert\letsymb{s},t_0)=(q_{jk})_{d\times\letsymb{d}}; \; q_{jk}=p(x_{j},t\vert\letsymb{x}_{k},s;X_{t_{0}}=x_{l}),\] 
which represent transition probabilities from a state at time $s$ to a state at a later time $t$, conditioned on the fixed initial state $X_{t_0} = x_l$. The subscript and superscript in $Q_{1}^{(l)}(t \vert \letsymb{s}, t_0)$ denote, respectively, the memory order (one-point) and the initial state on which the transition is conditioned.
Accordingly, Eq.~(\ref{joint_prob_order_3_eq1}) can be written as
\begin{equation}\label{joint_prob_order_3_eq2}
	p(x_j,t;x_k,s;x_{l},t_{0})=Q^{(l)}_{1}(t\vert\letsymb{s,t_{0}})_{jk}\Lambda(s\vert\letsymb{t_{0}})_{kl}q_{l}. 
\end{equation}

\par Summing both sides of the Eq.~(\ref{joint_prob_order_3_eq2}) over the initial probability and dividing it by the intermediate probability leads to 
\begin{equation}\label{intermeditate_cond_prob_eq3}
	p(x_j,t\vert\letsymb{x_k},s)=\sum_{l=1}^{d}q_{l}Q^{(l)}_{1}(t\vert\letsymb{s,t_{0}})_{jk}\frac{\Lambda(s\vert\letsymb{t_{0}})_{kl}}{p(x_{k},s)}
\end{equation}
On using the expression of the intermediate transition matrix $\Lambda(t\vert\letsymb{s})_{jk}=p(x_j,t\vert\letsymb{x_k},s)$, we finally obtain from Eq.~(\ref{intermeditate_cond_prob_eq3}) the compact form
\begin{equation}\label{eq04:intermediate transition matrix}
	\Lambda(t\vert\letsymb{s})=\sum_{l=1}Q^{(l)}_{1}(t\vert\letsymb{s},t_{0})D^{(l)}(s\vert\letsymb{t_{0}})
\end{equation}
so that we can express the probability vector as $\mathbf{p}(t)=\Lambda(t\vert\letsymb{s})\mathbf{p}(s)$, where, noting that $p(x_{k},s)=\sum_{l=1}^{d}q_{l}\Lambda(s\vert\letsymb{t_0})_{kl}$, the elements of the diagonal matrices $D^{(l)}(s\vert\letsymb{t_{0}})=diag(a_1^{(l)},a_2^{(l)},\ldots,a_d^{(l)})$ have the following expression
\begin{equation}\label{diagonal_matrices}
	a_{j}^{(l)}=\frac{q_{l}\Lambda(s\vert\letsymb{t_0})_{jl}}{\sum_{k=1}^{d}q_{k}\Lambda(s\vert \letsymb{t_{0}})_{jk}}.
\end{equation}
The elements of $D^{(l)}(s\vert\letsymb{t_0})$ can equivalently be written as $a_{j}^{(l)} = p(x_l,t_0 \vert \letsymb{x_j},s)$, which shows that Eq.~(\ref{eq04:intermediate transition matrix}) follows from an application of Bayes’ rule to the conditional probabilities.
\par Markov processes are characterized by the following Markov condition on the conditional probabilities
\begin{eqnarray*}
	p(x_{j_m},t_m\vert\letsymb{x_{j_{m-1}},t_{m-1}};\ldots;x_{j_1},t_{1}; x_{j_0},t_{0})\\
	=p(x_{j_m},t_m\vert\letsymb{x_{j_{m-1}},t_{m-1}}),
\end{eqnarray*}
which leads to the celebrated $\mathrm{CK}$-equation
\begin{equation}\label{Chapman-Kolmogorov_eq}
	p(x,t\vert\letsymb{z},s)=\sum_{y\in\mathrm{E}}p(x,t\vert\letsymb{y},t^{\prime})p(y,t^{\prime}\vert\letsymb{z},s)
\end{equation}
for all $t_0\leq\letsymb{s}\leq\letsymb{t^{\prime}}\leq\letsymb{t}$. Eq.~(\ref{Chapman-Kolmogorov_eq}) can be written in matrix form as $\Lambda(t\vert\letsymb{s})=\Lambda(t\vert\letsymb{t^{\prime}})\Lambda(t^{\prime}\vert\letsymb{s})$.
\par Accordingly, the one-point memory transition matrices should  satisfy  for all $t_{0}\leq\letsymb{s}\leq\letsymb{t}$ the following
\begin{description} 
	 \item[Consistency Condition] \begin{sloppypar}
		If the stochastic process $X_{t}$ is Markovian, then  $Q^{(1)}(t\vert\letsymb{s},t_{0})=Q^{(2)}(t\vert\letsymb{s},t_{0})=\cdots=Q^{(d)}(t\vert\letsymb{s},t_{0})=\Lambda(t\vert\letsymb{s})$ (see~Eq.~(27) in Ref.~\cite{CanturkBreuer2024}).
	\end{sloppypar} 
\end{description}
\par It is evident that $\sum_{l=1}^{d}D^{(l)}(s\vert\letsymb{t_{0}})=I$  which simplifies the intermediate transition matrix given by Eq.~(\ref{eq04:intermediate transition matrix}) when the one-point transition matrices are the same; that is, when the process is Markovian. One can then observe that invoking the  consistency condition in Eq.~(\ref{eq04:intermediate transition matrix}) indeed eliminates the effect of the initial state from the intermediate transition matrix. 
\par  We point out that $\lim_{s\rightarrow \letsymb{t_0}}D^{(l)}(s\vert t_{0})=diag(0,\ldots,1,\ldots,0)$ where the entry equal to $1$ appears in the $l^{th}$ diagonal position. Bearing in mind this fact and employing  $	\lim_{s\rightarrow\letsymb{t_{0}}}\Lambda(t\vert\letsymb{s})=\Lambda(t\vert t_{0})$, we obtain from Eq.~(\ref{eq04:intermediate transition matrix})  the  equality
\begin{equation}\label{transition_matrix_eq8}
Col_{j}(\Lambda(t\vert t_{0}))=Col_{j}(Q^{(j)}_{1}(t\vert s,t_{0})), \; j=1,2,\ldots,d,
\end{equation}
where $Col_{j}(A)$ denotes the $j^{th}$ column of a matrix $A$.
  Now, we assume that for all $t_{k}\in \mathrm{T}$ the conditional probabilities are dependent on the last and initial events:
 Then, the joint probability of order $m$ for the ordered time instants $t_{0}\leq\letsymb{t_{1}}\leq\ldots\letsymb{t_{m-1}}\leq\letsymb{t_{m}}$
 can be rewritten in terms of conditional probabilities as 
 \begin{eqnarray}\label{joint_prob_order_m_eq9}
 	&&\!\!\!\!\letsymb{p}_{m}(x_{j_m},t_m;x_{j_{m-1}},t_{m-1}; x_{j_{m-2}},t_{m-2};\ldots;x_{j_1},t_{1}; x_{j_0},t_{0})\nonumber\\
 	&&\quad =p(x_{j_m},t_m\vert\letsymb{x_{j_{m-1}},t_{m-1}}; x_{j_{0}},t_{0})\nonumber \\
 	&&\quad\times\letsymb{p}(x_{j_{m-1}},t_{m-1}\vert x_{j_{m-2}},t_{m-2}; x_{j_0},t_{0})\nonumber\\
 	&&\quad \times\; \ldots\; \times\;p(x_{j_1},t_{1}\vert \letsymb{x_{j_0}},t_{0})p(x_{j_0},t_{0})\nonumber\\
 	&&\quad = Q_{1}^{(j_0)}(t_m\vert t_{m-1},t_0)_{j_m,j_{m-1}}\nonumber\\
 	&&\quad \times\letsymb{Q}_{1}^{(j_0)}(t_{m-1}\vert t_{m-2},t_0)_{j_{m-1},j_{m-2}}\nonumber\\
 	&&\quad \times\ldots\times\letsymb{\Lambda}(t_1\vert t_0)_{j_1,j_0}p(x_{j_0},t_0)
 \end{eqnarray}
Invoking the Kolmogorov consistency condition
\begin{eqnarray*}
	&&\!\!\!\!\letsymb{p}_{m}(x_{j_m},t_m; x_{j_{m-2}},t_{m-2};\ldots;x_{j_1},t_{1}; x_{j_0},t_{0})\nonumber\\
	&& =\sum_{j_{m-1}}\letsymb{p}_{m}(x_{j_m},t_m;x_{j_{m-1}},t_{m-1}; \ldots;x_{j_1},t_{1}; x_{j_0},t_{0})\nonumber\\
\end{eqnarray*}
in Eq.~(\ref{joint_prob_order_m_eq9}) imposes  the condition, 
\begin{eqnarray}\label{condition_on_memory_transition_matrix}
	&&Q_{1}^{(j_0)}(t_{m}\vert t_{m-2},t_0)\nonumber\\
	&&\quad=	Q_{1}^{(j_0)}(t_{m}\vert t_{m-1},t_0)	Q_{1}^{(j_0)}(t_{m-1}\vert t_{m-2},t_0)
\end{eqnarray}
for all $t_{0}\leq\letsymb{t}_{m-2}\leq\letsymb{t_{m-1}}\leq\letsymb{t_m}$. This defines a $\mathrm{CK}$-equation type composition law for the one-point memory transition matrices. 
\par Once the one-point transition matrices are constructed consistently, the whole process is determined completely. Indeed, after constructing $\{Q^{(l)}(t\vert\letsymb{s},t_{0}),l=1,2,\ldots,d\}$, one can first  obtain the transition matrix $\Lambda(t\vert t_0)$ with the help of Eq.~(\ref{transition_matrix_eq8}), and then, the joint probability of order $m$ by means of Eq.~(\ref{joint_prob_order_m_eq9}). Moreover, knowledge of $\Lambda(t\vert t_0)$ allows one to determine the marginal probability vector $\mathbf{p}(t)=\Lambda(t\vert \letsymb{t_0})\mathbf{q}(t_0)$, and hence the moments of the process.
\par In the next section, we construct such a stochastic process on a two-dimensional sample space, evolving continuously in time.
\section{Construction of the process}\label{construct_the_process}
Let us consider the stochastic process $\{X_t\}_{t\in\mathrm{T}}$ with the sample space $\mathrm{E}=\{x_1,x_2\}$. According to Eq.~(\ref{condition_on_memory_transition_matrix}), for each fixed initial state, the process behaves similar to a Markov process. This observation enables us to construct each one-point memory transition matrix as the transition matrix of an associated Markov process. For a detailed treatment of the construction of transition matrices for Markov processes, the reader is referred to the classic works of Feller~\cite{Feller1970V2} and Gardiner~\cite{Gardiner2009}. To this end, we introduce the following time-homogeneous, one-point memory transition matrices:
\begin{subequations}\label{one_point_transition_matrices_eq11}
\begin{equation}
\!\!\!Q^{(1)}_{1}(t\vert\letsymb{s},t_{0})=\!\begin{bmatrix}
	k+(1-k)\e^{-\alpha\tau_{2}} & k(1-\e^{-\alpha\tau_{2}})\\
	(1-k)(1-\e^{-\alpha\tau_{2}}) & (1-k)+k\e^{-\alpha\tau_{2}}
\end{bmatrix}
\end{equation}
\begin{equation}
\!\!\!Q^{(2)}_{1}(t\vert\letsymb{s},t_{0})=\!\begin{bmatrix}
	(1-r)+r\e^{-\alpha\tau_{2}} & (1-r)(1-\e^{-\alpha\tau_{2}})\\
	r(1-\e^{-\alpha\tau_{2}}) & r+(1-r)\e^{-\alpha\tau_{2}}
\end{bmatrix}
\end{equation}
\end{subequations}
for $t_{0}\leq \letsymb{s}\leq \letsymb{t}$, where $\tau_{2}=t-s$ and $0\leq \letsymb{k,r}\leq 1$. The non-negative real parameter $\alpha$ characterizes the relaxation rate of the process, while the parameters $k$ and $r$ determine the transition rates between the states. We note that each $Q^{(l)}_{1}(t\vert \letsymb{s},t_{0})$ satisfies Eq.~(\ref{condition_on_memory_transition_matrix}). On using Eq.~(\ref{transition_matrix_eq8}) we find 
\begin{equation}
	\!\!\!\Lambda(t\vert\letsymb{t_{0}})=\begin{bmatrix}
		k+(1-k)\e^{-\alpha\tau_{1}} & (1-r)(1-\e^{-\alpha\tau_{1}})\\
		(1-k)(1-\e^{-\alpha\tau_{1}}) &r+(1-r)\e^{-\alpha\tau_{1}})
	\end{bmatrix},
\end{equation}  
where $\tau_{1}=t-t_0$.
\par  As noted above, knowledge of the matrices $\{Q^{(1)}_{1}(t\vert\letsymb{s},t_{0}),Q^{(2)}_{1}(t\vert\letsymb{s},t_{0}),\Lambda(t\vert\letsymb{t_{0}})\}$ is sufficient to specify the process. It is evident that the constructed process is irreducible and becomes Markovian iff $k+r=1$. Note that $det(\Lambda(t\vert\letsymb{t_0}))=k+r-1+(2-k-r)\e^{-\alpha(t-t_0)}$, which can be zero only if $k+r<1$. From now on, we assume that $k+r\geq\letsymb{1}$.  In order to proceed, the following remark is in order~\cite{Vacchini-Smirne-Breuer2011}:
\begin{remark}
	Let  $\{X_{t}\}_{t\in\mathrm{T}}$ be a stochastic process and  $\Lambda(t\vert\letsymb{t}_{0})$ the associated stochastic matrix on a finite samples space. The process is called 
	\textit{positively divisible (classical P-divisible)} if $\Lambda(t\vert\letsymb{t}_{0})=\tilde{\Lambda}(t\vert\letsymb{s})\Lambda(s\vert\letsymb{t}_{0})$ such that 
	$\tilde{\Lambda}(t\vert\letsymb{s})$ is also a stochastic matrix for all $t_0\leq s\leq t$.
\end{remark}
Classical P-divisibility is the generalization of $\mathrm{CK}$-equation to stochastic matrices. Since $\Lambda(t\vert\letsymb{t_0})$ is invertible, we can introduce the matrix $\Gamma(t\vert\letsymb{s})=\Lambda(t\vert\letsymb{t_0})\Lambda^{-1}(s\vert\letsymb{t_0})$, whose elements are
\begin{subequations}
	\begin{equation}
		\Gamma(t\vert\letsymb{s})_{11}=1-\frac{(1-k)\e^{-\alpha(s-t_0)}(1-\e^{-\alpha(t-s)})}{det(\Lambda(s\vert\letsymb{t_0}))},
	\end{equation}
	\begin{equation}
		\Gamma(t\vert\letsymb{s})_{22}=1-\frac{(1-r)\e^{-\alpha(s-t_0)}(1-\e^{-\alpha(t-s)})}{det(\Lambda(s\vert\letsymb{t_0}))},
	\end{equation}
	\begin{eqnarray}
		\Gamma(t\vert\letsymb{s})_{21}=1-\Gamma(t\vert\letsymb{s})_{11}, \;  \Gamma(t\vert\letsymb{s})_{12}=1-\Gamma(t\vert\letsymb{s})_{22}.
	\end{eqnarray}
\end{subequations}
Note that $\Gamma(t\vert \letsymb{s})$ is a stochastic matrix for all $t_0 \leq \letsymb{s} \leq \letsymb{t}$. It follows that $\Lambda(t\vert\letsymb{t_0}) = \Gamma(t\vert\letsymb{s}) \Lambda(s\vert\letsymb{t_0})$, implying that the constructed process is \textit{classical P-divisible}. We emphasize, however, that $\Lambda(t\vert\letsymb{s})$ in Eq.~\eqref{eq04:intermediate transition matrix} does not coincide with $\Gamma(t\vert\letsymb{s})$ unless $k + r = 1$, i.e., the condition under which the process becomes Markovian. In other words,  $\Lambda(t\vert\letsymb{t_0}) \neq \Lambda(t\vert\letsymb{s}) \Lambda(s\vert\letsymb{t_0})$. This observation shows that classical P-divisibility alone is not sufficient to determine whether a process is Markovian as also noted in Refs.~\cite{HanggiThomas1977,Vacchini-Smirne-Breuer2011}. One might argue that $\Lambda(t\vert\letsymb{s})$ could be replaced by $\Gamma(t\vert\letsymb{s})$ without affecting the statistical analysis. However, as we demonstrate in Sec.~\ref{statistical_anlaysis}, this is not the case: memory of the initial state manifests itself significantly in the statistical variables.
\par Using the treatment 
\begin{eqnarray*}
	\frac{d\mathbf{p}(t)}{d\letsymb{t}}&&=\lim_{\Delta\letsymb{t}\rightarrow\letsymb{0}}\frac{\mathbf{p}(t+\Delta\letsymb{t})-\mathbf{p}(t)}{\Delta\letsymb{t}}\\
	&&=\lim_{\Delta\letsymb{t}\rightarrow\letsymb{0}}\frac{\Lambda(t+\Delta\letsymb{t}\vert\letsymb{t_0})\mathbf{q}(t_0)-\Lambda(t\vert\letsymb{t_0})\mathbf{q}(t_0)}{\Delta\letsymb{t}}\\
	&&=\lim_{\Delta\letsymb{t}\rightarrow\letsymb{0}}\frac{\left(\Gamma(t+\Delta\letsymb{t}\vert\letsymb{
		t})-I\right)\mathbf{p}(t)}{\Delta\letsymb{t}}
\end{eqnarray*}
and defining $\mathcal{L}(t):=\lim_{\Delta\letsymb{t}\rightarrow\letsymb{0}}\frac{1}{\Delta\letsymb{t}}(\Gamma(t+\Delta\letsymb{t}\vert\letsymb{t})-I)$, we find the following master equation 
\begin{equation}\label{master_equation_eq13}
	\frac{d\mathbf{p}(t)}{d\letsymb{t}}=\mathcal{L}(t,t_0)\mathbf{p}(t).
\end{equation}
The explicit form of $\mathcal{L}(t,t_0)$ is  
\begin{equation}
	\mathcal{L}(t,t_0)=w(t,t_0)\mathrm{L}=w(t,t_0)\begin{bmatrix}
		-(1-k) & 1-r\\
		1-k & -(1-r)
	\end{bmatrix},
\end{equation}
where $w(t,t_0)=\frac{\alpha\e^{-\alpha(t-t_0)}}{det(\Lambda(t\vert\letsymb{t_0}))}$, which is non-negative for all $t\geq\letsymb{t_0}$. Note that $w(t,t_0)=\alpha$ when $k+r=1$, i.e., when the process is Markovian. Applying the shifting $t\rightarrow\tau=t-t_0$,  Eq.~(\ref{master_equation_eq13}) becomes  a time-local master equation with respect to $\tau$, since $w(t,t_0)\rightarrow w(\tau)$ and the generator reduces to $\mathcal{L}(\tau)$.     Defining the transition rates $W_{21}(\tau)=(1-k)w(\tau)$, $W_{12}(\tau)=(1-r)w(\tau)$, Eq.~(\ref{master_equation_eq13}) can be equivalently written as 
\begin{equation}\label{master_eq_marginal_prob}
	\frac{d\letsymb{p}(x_1,\tau)}{d\letsymb{\tau}}=W_{12}(\tau)\letsymb{p}(x_2,\tau)-W_{21}(\tau)\letsymb{p}(x_1,\tau),
\end{equation}
which demonstrates the fundamental fact that there are non-Markovian processes  satisfying a master equation (or differential $\mathrm{CK}$-equation) with non-negative, time-dependent transition rates. It is important to note that time-inhomogeneous (i.e., non-homogeneous) Markov processes also obey a differential  $\mathrm{CK}$-equation with non-negative, time-dependent transition rates~\cite{Vassiliou2023}.
\par  We refer to the process constructed above as the Funessian process, inspired by Jorge Luis Borges’s story \textit{Funes el memorioso}~\cite{Borges1962}. In Borges’s story, Ireneo Funes develops a perfect memory after a traumatic fall from his horse. As shown in Section~\ref{statistical_anlaysis}, the memory of the first event—which plays the role of Funes’s accident—fundamentally alters the nature of the process.

\section{Statistical analysis of the process}\label{statistical_anlaysis}
We calculate the mean value, variance, correlation function of the process, and mutual information at stationary state. To this end, we follow their formulations introduced by Gardiner (see~Ref.~\cite{Gardiner2009}, pp.~75-76). The stationary value of the mean value conditioned on the initial state $x_{0}\in\mathrm{E}$
\[\braket{X_{t}}{[x_0,t_0]}=\sum_{j=1}^{2}x_jp(x_j,t\vert\letsymb{x_0,t_0})\]
is defined as $\braket{X_t}{\letsymb{x_0}}_{st}:=\lim_{t_0\rightarrow\letsymb{-\infty}}\braket{X_{t}}{[x_0,t_0]}$. Then, a straightforward calculation yields
\begin{eqnarray}
	\braket{X_t}{\letsymb{x_0}}_{st}&&= (k\delta_{x_1,x_0}+(1-r)\delta_{x_2,x_0})x_1\nonumber\\
	&&+((1-k)\delta_{x_1,x_0}+r\delta_{x_2,x_0})x_2,
\end{eqnarray}
where $\delta_{x_j,x_0}$ is Kronecker delta. Similarly, the stationary value of the second moment conditioned on the initial state
\[\braket{X_{t}^{2}}{[x_0,t_0]}=\sum_{j=1}^{2}x_{j}^{2}p(x_j,t\vert\letsymb{x_0,t_0})\]
is defined as $\braket{X_t^{2}}{\letsymb{x_0}}_{st}:=\lim_{t_0\rightarrow\letsymb{-\infty}}\braket{X_{t}^{2}}{[x_0,t_0]}$. Then, the stationary value of the variance is
\begin{eqnarray}
&&\!\!\!\!Var(X_t\vert\letsymb{x_0})_{st}=	\braket{X_t^{2}}{\letsymb{x_0}}_{st}-\braket{X_t}{\letsymb{x_0}}_{st}^{2}\nonumber\\
&& = (x_{1}-x_{2})^{2}(k(1-k)\delta_{x_{1},x_{0}}+r(1-r)\delta_{x_{2},x_{0}}).
\end{eqnarray}
The dependence of the mean value and variance on the initial state is evident. Notably, while the variance loses its dependence on the initial state when $k = r$, the mean value retains such dependence for all values of $k$ and $r$ satisfying $k + r > 1$. When $k+r=1$, corresponding to a Markov process, these quantities reduce to
\begin{eqnarray*}
	\braket{X_t}{\letsymb{x_0}}_{st}=kx_1+rx_2; \;Var(X_t\vert\letsymb{x_0})_{st}=kr(x_1-x_2)^{2}
\end{eqnarray*}
 which are independent of the initial state, as expected.   From the definition of two-time correlation,
\begin{eqnarray}\label{correlation_function_calculation}
		&&\braket{X_{t}X_{s}}{[x_{j_{0}},t_{0}]}=\sum_{j,k=1}^{2}x_{j}x_{k}p(x_{j},t;x_{k},s\vert\letsymb{x_{j_{0}},t_{0}})\nonumber\\
		&&\quad\quad =\sum_{j,k=1}^{2}x_{j}x_{k}p(x_{j},t\vert\letsymb{x}_{k},s;x_{j_0},t_0)p(x_{k},s\vert\letsymb{x}_{j_{0}},t_{0})\nonumber\\
		&&\quad\quad =\sum_{j,k=1}^{2}x_{j}x_{k}Q_{1}^{(j_0)}(t\vert\letsymb{s},t_0)_{jk}\Lambda(s\vert\letsymb{t_0})_{kj_{0}}, 
\end{eqnarray}
the stationary value of the correlation function 
\[	C(t,s\vert\letsymb{x_{j_{0}}})_{st}=\lim_{t_{0}\rightarrow\letsymb{-\infty}}\left(\braket{X_{t}X_{s}}{[x_{j_0},t_{0}]}-\braket{X_{t}^{2}}{[x_{j_0},t_0]}\right)\]
 can be obtained as 
\begin{equation}\label{NonMarkov_correlation_function}
\begin{aligned}
	C(t,s\vert\letsymb{x_{j_0}})_{st}&=(x_{1}-x_{2})^{2}(k(1-k)\delta_{x_{1},x_{j_0}}\\
		&+r(1-r)\delta_{x_{2},x_{j_0}})\e^{-\alpha(t-s)}.
\end{aligned}
\end{equation}
Taking average over the initial state yields
\begin{equation}\label{correlation_function}
\begin{aligned}
&\tilde{C}(t,s)_{st}=\sum_{j_{0}=1}^{2}q_{j_{0}}C(t,s\vert\letsymb{x_{j_0}})_{st}\\
	&\quad =(x_{1}-x_{2})^{2}(k(1-k)q_{1}+r(1-r)q_{2})\e^{-\alpha(t-s)}.
\end{aligned}
\end{equation}
If we replaced  $Q^{(j_{0})}_{1}(t\vert\letsymb{s},t_0)_{jk}$ with $\Gamma(t\vert\letsymb{s})_{jk}$  in the calculation of the correlation function in Eq.~(\ref{correlation_function_calculation}), we would find
\[C(t,s\vert\letsymb{x_{j_0}})_{st}=(k(1-k)\delta_{x_1,x_{j_0}}+r(1-r)\delta_{x_2,x_{j_0}})(x_1-x_2)^{2},\]
which is equal to $Var(X_t\vert\letsymb{x_0})_{st}$ and constant with respect to time difference $\tau=t-s$, implying that the process would be perfectly correlated in time. This result is, however, fundamentally different from that  given by Eq.~(\ref{NonMarkov_correlation_function}).  This implies that replacing $\Lambda(t\vert\letsymb{s})$ with $\Gamma(t\vert\letsymb{s})$ captures neither the true statistical properties of the process nor its physical structure.
\par
\begin{figure}[h]
	\centering
	\begin{subfigure}[b]{0.50\columnwidth}
		\centering
		\includegraphics[width=\textwidth]{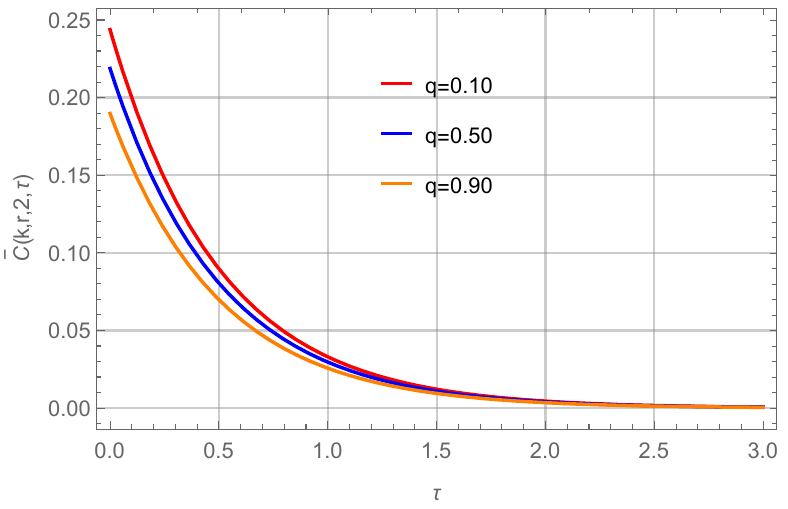}
		\caption{}
		\label{fig_corr_func_differ_kr}
	\end{subfigure}%
	\hfill
	\begin{subfigure}[b]{0.50\columnwidth}
		\centering
		\includegraphics[width=\textwidth]{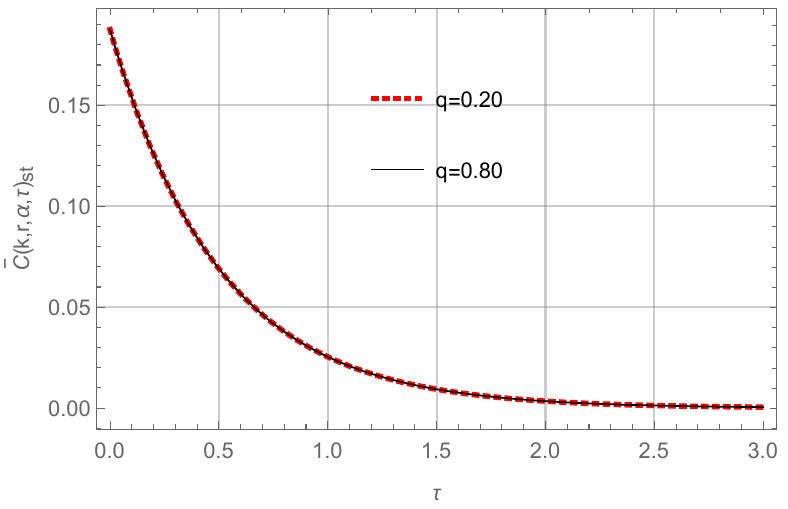}
		\caption{}
		\label{fig_corr_func_same_kr}
	\end{subfigure}
	\captionsetup{font=footnotesize,textfont=it}
	\caption{Average correlation function vs.\ time difference $\tau$, for different initial distribution $\mathbf{q}(t_0)=(q,1-q)$ (time unit set to 1 and $\alpha=2$). (a)  $k=0.75, r=0.50$; (b)  $k=r=0.75$.}
	\label{fig_corr_func}
\end{figure}

  When the process becomes Markovian (i.e., k+r=1), the correlation function reduces to 
\begin{equation}
	\tilde{C}(t,s)_{st}=C(t,s\vert\letsymb{x_{j_0}})_{st}=kr(x_1-x_2)^{2}\e^{-\alpha(t-s)},
\end{equation} 
which no longer depends  on the initial state. The correlation function of the constructed process decays exponentially, as is characteristic of Markov processes. On the other hand, the non-Markovian nature of the process can still be revealed through the correlation function by altering the initial density of states. This fact can be captured in Figure~\ref{fig_corr_func_differ_kr} which illustrates the change of the average correlation function, $\tilde{C}(k,r,\alpha,\tau)_{st}=\tilde{C}(t,s)_{st}$, with respect to time difference $\tau=t-s$ when $k\neq\letsymb{r}$. 
However, when $k = r$, the correlation function fails to reveal the non-Markovian nature of the process, as illustrated in Fig.~\ref{fig_corr_func_same_kr}. In the following, we show that, even in this case, mutual information can uncover the underlying non-Markovian behavior. To this end, let us recall that  the conditional probability of a Markov chain $Z_{t_0}\rightarrow\letsymb{Y_{s}}\rightarrow\letsymb{X_{t}}$ has the property 
\begin{equation}\label{conditional_equality_Markovian}
	p(x,t;y,s\vert\letsymb{z},t_0)=p(x,t\vert\letsymb{y},s)p(y,s\vert\letsymb{z},t_0)
\end{equation}
for all $t_0\leq\letsymb{s}\leq\letsymb{t}$. Note that Eq.~(\ref{conditional_equality_Markovian}) fails for non-Markovian process, in particular, for the process constructed above. This fact can be expressed in terms of the difference of the corresponding Shannon entropies, averaged over the conditional data $z$, 
\begin{equation}\label{entropic_expression}
	\begin{aligned}
		&\sum_{z}p(z,t_0)\biggl(-\sum_{x,y}	p(x,t;y,s\vert\letsymb{z},t_0)\ln{	p(x,t;y,s\vert\letsymb{z},t_0)}\\
		&+\sum_{x,y}p(x,t\vert\letsymb{y},s)p(y,s\vert\letsymb{z},t_0)\ln(p(x,t\vert\letsymb{y},s)p(y,s\vert\letsymb{z},t_0))\biggr),
	\end{aligned}
\end{equation}
which is identically zero when the process is Markovian. After some algebra, Eq.~(\ref{entropic_expression}) reduces to an expression in terms of mutual information (see ~\cite{CoverThomas2006}, pp.~19-35) as

\begin{equation}\label{mutual_info_middle_expression}
	I(X_t;Y_s\vert\letsymb{Z_{t_0}})+I(X_t;Z_{t_0})-I(X_t;Y_s).
\end{equation}
 Using the chain rule $I(X_t;Y_s,Z_{t_0})=I(X_t;Y_s\vert\letsymb{Z}_{t_0})+I(X_t;\letsymb{Z_{t_0}})$ in Eq.~(\ref{mutual_info_middle_expression}), we finally obtain
\begin{equation}\label{mutual_infor_final_expression}
	\begin{aligned}
		I(X_t;Y_s,Z_{t_0})-I(X_t;Y_s)=I(X_{t};Z_{t_0}\vert\letsymb{Y_{s}}).
	\end{aligned}
\end{equation}
 Thus, the quantity $I(X_{t};Z_{t_0}\vert Y_{s})$ in Eq.~(\ref{mutual_infor_final_expression}) can be used to measure the non-Markovianity of the process; it simply measures the possible effect of the past information ($Z_{t_0}$) on the determination of the future sate ($X_{t}$) provided that we know the present state ($Y_{s}$). In other words, it is a quantitative answer to the following question: Once  $Y_{s}$ (the present) is known, does $Z_{t_0}$ (the past) give any extra information about $X_{t}$ (the future)? The answer is simply negative for a Markov process. We note that a similar approach, based on probability distributions associated with measurement outcomes, was proposed in Ref.~\cite{Budini2018} for quantum processes. Below we show that for the constructed process, the knowledge of the initial state (i.e., the past) can always give extra information about the state at a later time $t$ (i.e., the future) after a measurement at the intermediate time $s$ (i.e., the present), regardless of how much time has passed over the initial time. That is why we call the process Funessian, since the initial state affects the whole future. From now on, we replace $Y_{s}$ with $X_s$ and $Z_{t_0}$ with $X_{t_0}$. Let $\Lambda=(\lambda_{jk})_{d\times\letsymb{d}}$ be a stochastic matrix and $\mathbf{r}=(r_1,r_2,\ldots,r_d)$ be a probability vector. Then, introducing
\[H(\Lambda\vert\letsymb{\mathbf{r}}):=- \sum_{j,k=1}^{d}r_{k}\Lambda_{jk}\ln{\Lambda_{jk}},\]
the stationary value of mutual information  is obtained as
\begin{equation}\label{mutual_information_nonMarkovian}
	\begin{aligned}
		&\!\!I_{st}(X_{t}:X_{t_{0}}\vert\letsymb{X_{s}})=\lim_{t_{0}\rightarrow\letsymb{-\infty}}I(X_{t}:X_{t_{0}}\vert\letsymb{X_{s}})\\
		&= \sum_{l=1}^{2}q_{l}\left(H(\Lambda_{st}(t\vert\letsymb{s})\vert\mathbf{p}_{st})-H(Q_{1}^{(l)}(t\vert\letsymb{s,t_0})\vert\letsymb{\mathbf{r}^{(l)}})\right),
	\end{aligned}
\end{equation}
where, noting that $\Lambda_{st}(\cdot\vert\cdot)=\lim_{t_{0}\rightarrow\letsymb{-\infty}}\Lambda(\cdot\vert\cdot)$,  $\mathbf{r}^{(l)}=Col_{l}(\Lambda_{st}(s\vert_{t_0}))$ and $\mathbf{p}_{st}=\Lambda_{st}(s\vert\letsymb{t_0})\mathbf{q}({t_0})$.  In particular, $\mathbf{r}^{(1)}=(k,1-k)$ and $\mathbf{r}^{(2)}=(1-r,r)$.  One needs the stationary values of the matrices $D^{(l)}(t\vert\letsymb{t_0})$ to find $\Lambda_{st}(t\vert\letsymb{s})$ in Eq.~(\ref{eq04:intermediate transition matrix}).   $D^{(l)}_{st}(t\vert\letsymb{t_0})$ can be obtained from Eq.~(\ref{diagonal_matrices}).  
\begin{figure}[h]
	\centering
	\begin{subfigure}[b]{0.50\columnwidth}
		\centering
		\includegraphics[width=\textwidth]{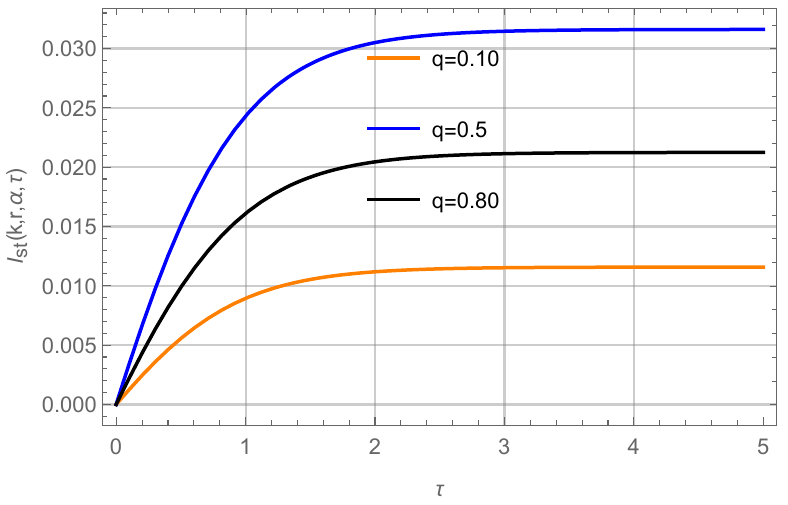}
		\caption{}
		\label{fig_mutual_info_differ_kr}
	\end{subfigure}%
	\hfill
	\begin{subfigure}[b]{0.50\columnwidth}
		\centering
		\includegraphics[width=\textwidth]{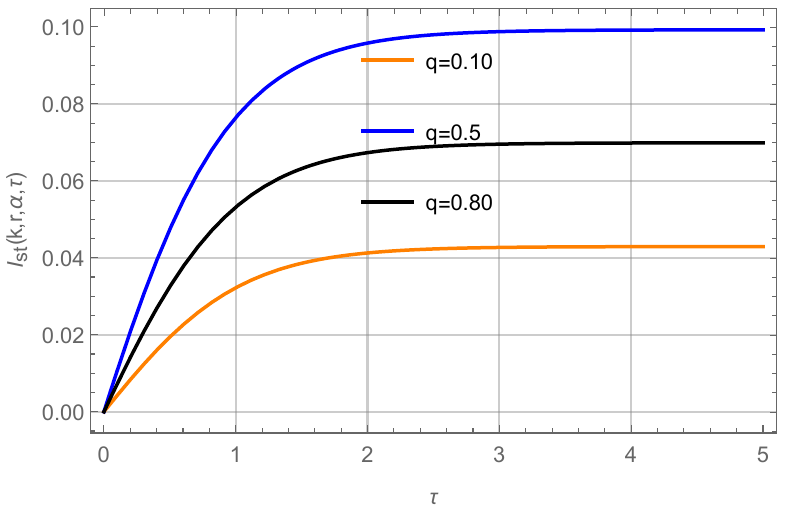}
		\caption{}
		\label{fig_mutual_info_same_kr}
	\end{subfigure}
	\captionsetup{font=footnotesize,textfont=it}
	\caption{Mutual information vs.\ time difference $\tau$, for different initial distribution $\mathbf{q}(t_0)=(q,1-q)$ (time unit set to 1 and $\alpha=2$). (a) $ k=0.75, r=0.50$; (b) $k=r=0.75$.}
	\label{fig_mutual_info_non_Markovian}
\end{figure}
  \par
\begin{figure}[h]
   	\centering
   	\begin{subfigure}[b]{0.50\columnwidth}
   		\centering
   		\includegraphics[width=\textwidth]{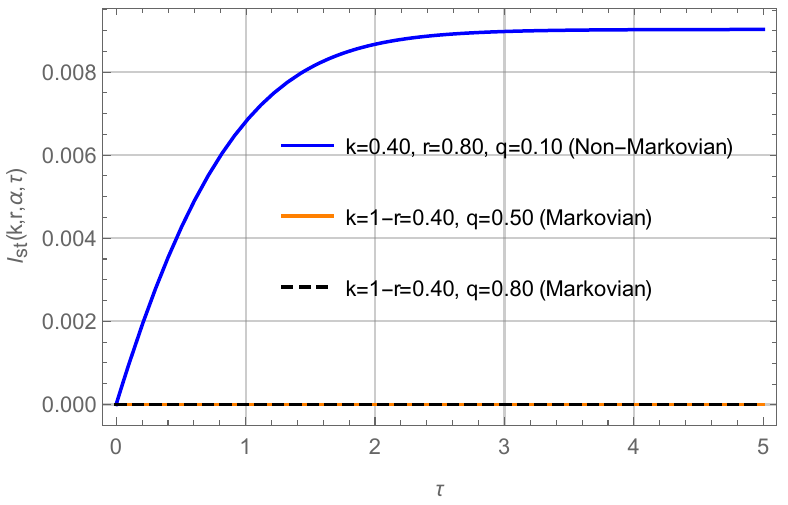}
   		\caption{}
   		\label{fig_mutual_info_Markov_nonMarkov}
   	\end{subfigure}%
   	\hfill
   	\begin{subfigure}[b]{0.50\columnwidth}
   		\centering
   		\includegraphics[width=\textwidth]{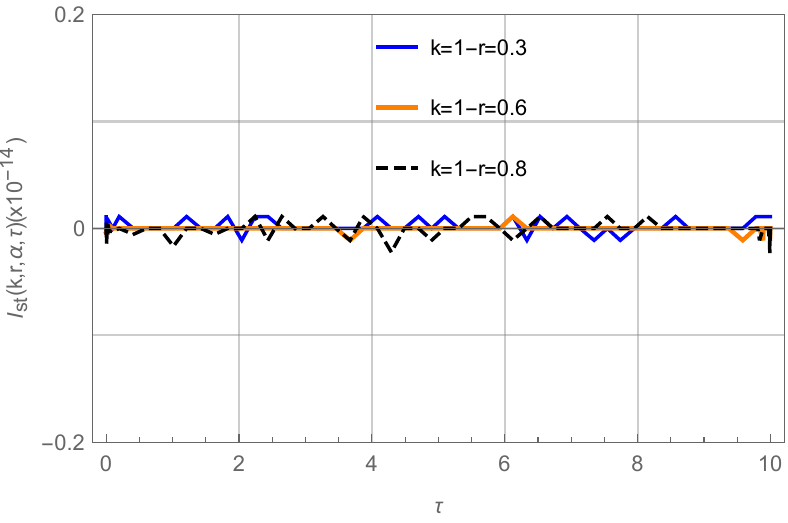}
   		\caption{}
   		\label{fig_mutual_info_Markov}
   	\end{subfigure}
   	\captionsetup{font=footnotesize,textfont=it}
   	\caption{Mutual information vs.\ time difference $\tau$, for different initial distribution $\mathbf{q}(t_0)=(q,1-q)$ (time unit set to 1 and $\alpha=2$). (a) Non-Markovian vs.\ Markovian case; (b) Markovian case for different $k$ and fixed $q=0.6$.}
   	\label{fig_compare_mutual_info_Markov_nonMarkov}
   \end{figure}

  We note that  $ I_{st}(X_{t};X_{t_{0}}\vert\letsymb{X_{s}})$ is a function of time difference $\tau=t-s$. When $\tau = 0$, the future coincides with the known present. Consequently, knowledge of the past provides no additional information about the future, since the future state is already known. Therefore, $I_{st}(X_{s};X_{t_{0}}\vert\letsymb{X_{s}})=0$, regardless of whether or not the process is Markovian. This is clearly observed  in Figures~\ref{fig_mutual_info_non_Markovian} and~\ref{fig_compare_mutual_info_Markov_nonMarkov}, which illustrate the change of mutual information, $I_{st}(k,r,\alpha,\tau)=I_{st}(X_{t};X_{t_{0}}\vert\letsymb{X_{s}})$, as a function of time difference $\tau$. In contrast to correlation function, Figure~\ref{fig_mutual_info_non_Markovian} shows that the non-Markovian behavior of the process is significantly detected by the mutual information through changes in the initial probability, regardless of whether $k$ and $r$ are equal.
  \par For Markovian case (i.e., k+r=1),  $\mathbf{r}^{(1)}=\mathbf{r}^{(2)}=\mathbf{p}_{st}=(k,1-k)$ and $Q_{1}^{(1)}(t\vert\letsymb{s},t_0)=Q_{2}^{(1)}(t\vert\letsymb{s},t_0)=\Lambda_{st}(t\vert\letsymb{s})$ so that $I_{st}(X_{t};X_{t_{0}}\vert\letsymb{X_{s}})=0$. This result is illustrated in both Figures~\ref{fig_mutual_info_Markov_nonMarkov} and~\ref{fig_mutual_info_Markov}. The fluctuations in mutual information for Markovian case in Figure~\ref{fig_mutual_info_Markov} are  of order ${10^{-14}}$ which is practically zero.
\par   Figure~\ref{fig_mutual_info_non_Markovian} indicates that memory of the initial state in the stochastic process persists over time. This suggests that non-Markovian processes cannot be treated as perturbations of Markov processes, in agreement with the argument of Van Kampen~\cite{VanKampen1998}. On the other hand, they may still exhibit similarities to Markov processes; for example, their correlation functions can behave similarly to those of Markov processes.  
\section{Application to Random Walk}\label{application_random_walk}
We consider the random walk
\begin{equation}\label{random_walk_form}
	S(t)=X_{t_0}+\sum_{k=1}^{N(t)}X_{t_k}.
\end{equation}
constructed from the stochastic process $X_{t}$ presented above, where $N(t)$ is a Poisson process with rate $\lambda$ governing the jumps and $X_t\in\{x_1,x_2\}$ is evaluated at jump times $\{t_k\}$. In the orthodox continuous-time random walk (CTRW) approach, the statistics of $S(t)$ is derived from a joint distribution $\psi(x,t)$ which is typically factorized into two independent probability density functions (pdf), $\psi(x,t)=w(t)g(x)$, where $w(t)$ and $g(x)$ denote the waiting-time and jump-length distributions, respectively \cite{MontrollWeiss1965, BouchaudGeorges1990}. Equation~(\ref{random_walk_form}) instead assumes a correlation between the jump length and the jumping time, since $X_{t_m}$ and $X_{t_n}$ are correlated through the transition matrix $\Lambda(t_m\vert\letsymb{t_n})$, where $t_m$ and $t_n$ are jumping times. Such correlated dynamics have also been proposed~\cite{MonteroMasoliver2007} and utilized as a remedy for the infinite-velocity problem arising in Lévy flights, leading to the concept of Lévy walks~\cite{Zaburdaev2006,ZaburdaevChukbar2002,ZaburdaevDenisov2015}.    Furthermore,  non-Markovian generalizations of random walks are commonly formulated in a way that memory effects arise from non-exponential waiting-time distributions~\cite{AfekDavidson2023,ZaburdaevDenisov2015}. In contrast, in the present work the sampling mechanism remains Poissonian, while memory is introduced through the intrinsic dynamics of the underlying process $X_t$ which governs the jump lengths. This leads to a distinct class of stochastic processes in which correlations between increments originate from the temporal structure of $X_t$, rather than from the renewal properties of the random walk. We note that the random walk presented above is also different from Elephant Random Walk (ERW) which is self-interacting, thereby, memory effect of which is  internal and  depends on past trajectory, i.e., $S_{n+1}$ depends on $\{S_k\}_{k\leq\letsymb{n}}$~\cite{Roy2026,Gut2022}. However,  the memory effect of our model is external and dynamical, comes from $X_t$ with memory of $X_{t_0}$.
\par Let the sample space of $S(t)$ be $Z_{S}=\{z\vert z=nx_1+mx_2; n,m\in\mathbb{Z}\}$.  Over a short time interval $\left[t,t+\Delta\letsymb{t}\right]$, two events are possible: 1) one jump occurs with probability $\lambda\Delta\letsymb{t}$ and  $S(t)$ increases by a random amount $X_t$;  2) with probability $1-\lambda\Delta\letsymb{t}$ no jumps occurs and $S(t)$ remains unchanged. This leads to the balance equation 
 \[P(z,t+\Delta\letsymb{t})=(1-\lambda\Delta\letsymb{t})P(z,t)+\lambda\Delta\letsymb{t}\sum_{j=1}^{2}p(x_j,t)P(z-x_j,t).\]
 Taking the limit $\Delta\letsymb{t}\rightarrow\letsymb{0}$, we obtain the master equation 
 \begin{equation}
 	\frac{\partial}{\partial\letsymb{t}}P(z,t)=\lambda\sum_{j=1}^{2}p(x_j,t)P(z-x_j,t)-\lambda\letsymb{P(z,t)}.
 \end{equation}
 Multiplying the master equation by $z$ and summing over $z$, we obtain  an equation for the mean value $\average{S(t)}$
 \begin{equation}
 	\frac{d\average{S(t)}}{d\letsymb{t}}=\lambda\sum_{j=1}^{2}p(x_j,t)\sum_{z}zP(z-x_j,t)-\lambda\average{S(t)}.
 \end{equation} 
 In the gain term, shifting the variable $z\rightarrow\letsymb{z'+x_j}$ gives
 \[\sum_{z}zP(z-x_j,t)=\sum_{z^{\prime}}(z^{\prime}+x_j)P(z^{\prime},t)=\average{S(t)}+x_j\]
 which, by substitution, yields
 \begin{equation}\label{master_eq_mean_random_walk}
 	\frac{d\average{S(t)}}{d\letsymb{t}}=\lambda\average{X_t}.
 \end{equation}
 Integrating Eq.~(\ref{master_eq_mean_random_walk}) gives the mean value 
 \begin{equation}\label{mean_value_random_walk_first_step}
 \average{S(t)}=\average{X_{t_0}}+\lambda\int_{t_0}^{t}\average{X_s}ds,
 \end{equation}
 where, $\average{X_{t_0}}=x_1q_1+x_2q_2$. Using Eq.~(\ref{master_eq_marginal_prob}), one finds 
 \begin{equation}
 	\begin{aligned}
 		\average{X_{s}}&=M_{1}+ (k+r-1)\average{X_{t_0}}\\
 		&+((2-k-r)\average{X_{t_0}}-M_1)\e^{-\alpha(s-t_0)}
 	\end{aligned}
 \end{equation}  
 where $M_{1}=(1-r)x_1+(1-k)x_2$. Substituting  $\average{X_s}$ into the integral of Eq.~(\ref{mean_value_random_walk_first_step}), we finally obtain
 \begin{equation}
 	\begin{aligned}
 		&\average{S_{t}}=\average{X_{t_0}}+\lambda\left(M_{1}+ (k+r-1)\average{X_{t_0}}\right)(t-t_0)
 		\\
 		&+((2-k-r)\average{X_{t_0}}-M_1)(\lambda/\alpha)(1-\e^{-\alpha(t-t_0)})
 	\end{aligned}
 \end{equation} 
After a long time, $t\gg\letsymb{t_0}$, the second term becomes dominant and one can roughly read
\[\average{S_{t}}\propto \lambda\left[M_1+(k+r-1)\average{X_{t_0}}\right](t-t_0)\]  
which carries the effect of the initial state through the factor $(k+r-1)\average{X_{t_0}}$ that disappears for the Markovian case, i.e., when $k+r=1$.
 \par  Following the same procedure for the second moment as for the mean value leads to
 \[\frac{d\average{S^{2}(t)}}{d\letsymb{t}}=2\lambda\average{S(t)}\average{X_t}+\lambda\average{X^{2}_{t}}.\]
 Then, the variance $\sigma(S(t))=\average{S^{2}(t)}-\average{S(t)}^{2}$ satisfies
 \begin{equation*}
 		\frac{d\sigma(S(t))}{d\letsymb{t}}=\frac{d\average{S^{2}(t)}}{d\letsymb{t}}-2\average{S(t)}\frac{d\average{S(t)}}{d\letsymb{t}}=\lambda\average{X^{2}_t},
 \end{equation*}
which, by integration, gives
 \begin{equation}\label{variance_random_walk_first_step}
 	\sigma(S(t))=\sigma(X_{t_0})+\lambda\int_{t_0}^{t}\average{X_{s}^{2}}ds.
 \end{equation}
 Using  Eq.~(\ref{master_eq_marginal_prob}) again,  one finds
 \begin{equation}
 	\begin{aligned}
 	\average{X_{s}^{2}}&=M_{2}+ (k+r-1)\average{X_{t_0}^{2}}\\
 	&+((2-k-r)\average{X_{t_0}^{2}}-M_2)\e^{-\alpha(s-t_0)}	,
 	\end{aligned}
 \end{equation}
 where $M_{2}=(1-r)x_1^{2}+(1-k)x_{2}^{2}$. Substituting  $\average{X_s^{2}}$ into the integral of Eq.~(\ref{variance_random_walk_first_step}), we finally obtain
 \begin{equation}
 	\begin{aligned}
 		&\sigma(S_t) = \sigma(X_{t_0})
 		+ \lambda\left(M_2 + (k+r-1)\langle X_{t_0}^2\rangle\right)(t-t_0) \\
 		&\quad + \bigl((2-k-r)\langle X_{t_0}^2\rangle - M_2\bigr)
 		\frac{\lambda}{\alpha}\left(1 - e^{-\alpha(t-t_0)}\right)
 	\end{aligned}
 \end{equation}
 After a long time, $t\gg\letsymb{t_0}$, the second term becomes dominant and one can roughly read
 \[\sigma(S_{t})\propto \lambda\left[M_2+(k+r-1)\average{X_{t_0}^{2}}\right](t-t_0)\]  
which carries the effect of the initial state through the factor $(k+r-1)\average{X_{t_0}^{2}}$ that disappears for the Markovian case, i.e., when $k+r=1$. The  effective diffusion coefficient is $D_{eff}=\lim_{t\rightarrow\infty}(1/2)\left(\partial\sigma(S(t))/\partial\letsymb{t}\right)=(1/2)\lambda[M_2+(k+r-1)\average{X_{t_0}^{2}}]$. Although $\sigma(S(t))$ grows linearly in time, indicating normal diffusive scaling, the diffusion coefficient depends explicitly on the second moment of the initial distribution, revealing a persistent memory of the initial preparation encoded in the non-Markovian dynamics of $X_{t}$. In standard Markovian diffusion, $D_{eff}$ is uniquely determined by the stationary distribution of the underlying process $X_{t}$, independent of initial conditions. In the present case, however, the asymptotic state $\mathbf{p}_{st}=\lim_{t\rightarrow \infty}\Lambda(t\vert\letsymb{t_0})$ of $X_{t}$ retains an explicit dependence on the initial distribution, $\mathbf{p}_{st}=(kq_1+(1-r)q_2,(1-k)q_1+rq_2)$. This breakdown of uniqueness implies a failure of ergodicity at the level of transport: long-time diffusion properties such as $D_{eff}$ do not become independent of the initial preparation. This suggests that memory effects provide a mechanism for renormalizing diffusion through the memory of the initial state.
\section{Results and Discussion}\label{conclusion_section}
We have constructed a continuous-time, two-dimensional, classical P-divisible non-Markovian stochastic process with explicit dependence on the initial state that satisfies the differential $\mathrm{CK}$-equation with non-negative transition rates.  This demonstrates that non-Markovian processes may satisfy the differential $\mathrm{CK}$-equation, a property traditionally associated with Markov processes. The underlying framework is sufficiently general to construct any continuous-time, finite-dimensional stochastic process with memory exclusively of the  initial state. Although the correlation function exhibits exponential decay similar to that of Markov processes, we showed that mutual information reliably detects the non-Markovian character of the process. We note that the constructed process is irreducible which is a necessary property for ergodicity. However, its application to random walks further reveals that memory of the initial state can break ergodicity and modify transport properties such as the diffusion coefficient. 

\par
It has recently been shown that the time evolution of finite-dimensional stochastic processes can be uniquely represented by quantum processes~\cite{Canturk2026a}. In this representation, probability vectors and stochastic matrices are mapped to quantum states and quantum operations, respectively (see Chp.~2 of~\cite{Watrous2018} for a self-contained introduction to quantum processes). This correspondence supports the widely accepted idea that classical P-divisibility and the differential $\mathrm{CK}$-equation extend respectively to completely positive divisibility and the quantum master equation. However, the process constructed in the present work demonstrates that classical P‑divisibility alone is insufficient to characterize Markovianity. Consequently, characterizing quantum Markovianity solely in terms of CP divisibility may be questionable. A detailed discussion of the implications of these results for proposed criteria of quantum Markovianity lies beyond the scope of this paper and will be addressed elsewhere~\cite{CanturkBagci2026b}.

\begin{acknowledgments}
B. C. and G. B. B. gratefully acknowledge financial support from TUBITAK (Project No. : 124C348). 
\end{acknowledgments}


\bibliography{apssamp}

\end{document}